\documentstyle[12pt,amsfonts]{article}
\def\lra{{\longrightarrow }}

\def\dem{{\noindent{\bf Proof. }}}
\def\qed{{\hfill $\square$\\}}

\newcommand{\C}{{\Bbb C}}
\newcommand{\Z}{{\Bbb Z}}

\newcommand{\Q}{{\Bbb Q}}
\newcommand{\PS}{{\Bbb P}}

\newcommand{\cC}{{\cal C}}

\newcommand{\cO}{{\cal O}}

\newtheorem{thm}{Theorem}[section]
\newtheorem{prop}[thm]{Proposition}
\newtheorem{cor}[thm]{Corollary}
\newtheorem{lemma}[thm]{Lemma}
\newtheorem{definition}[thm]{Definition}

\begin{document}
\title{\bf  MOTIVES OF UNIRULED 3-FOLDS}
\author{Pedro Luis del Angel and  
Stefan M\"uller-Stach }

\date{version July 1996}
\maketitle

\begin{abstract} {J. Murre has conjectured that every smooth projective
variety $X$ of dimension $d$ admits a decomposition of the diagonal 
$\Delta=p_0+...+p_{2d} \in CH^d(X \times X) \otimes \Q$ such that 
the cycles $p_i$ are idempotent correspondences which lift the
K\"unneth components of the identity map in \'etale cohomology. If
this decomposition induces a filtration on the Chow groups of $X$ which is
independent of all choices, we call it a Murre decomposition. In this 
paper we propose candidates for such projectors on 3-folds by using 
fiber structures. Using Mori theory,
we prove that every smooth uniruled complex 3-fold 
admits a Murre decomposition. \\
{\bf Key Words:} Chow groups, cohomology groups, Chow motive, extremal ray. }
\end{abstract}
\parindent=0pt

\section{Introduction}

Let $F$ be a subfield of $\C$. We denote by $V(F)$ the category of smooth, 
projective varieties over $F$ with the usual morphisms. Let $CV(F)$ be the
category with the same underlying object, but where the morphisms are
replaced by correspondences of degree zero, i.e. for two irreducible varieties $X,Y$ we have
$Mor(X,Y):=CH^{\dim(X)}(X \times Y)$. If $f \in Mor(X,Y)$ we view it as a
homomorphism $f_\ast: CH^*(X) \to CH^*(Y)$, by defining
$f_*(W)=(pr_2)_* ( (W \times X) \cap f) $.
Given $X_1,X_2,X_3 \in V(k)$ the composition of correspondences 
$f \in Mor(X_1,X_2)$ and $g \in Mor(X_2,X_3)$ is defined by
$$ g \circ f = (pr_{13})_*\{( pr_{12})^* f \cap (pr_{23})^* g \}  $$
An element $p \in Mor(X,X)$ is called a {\bf projector} if $p \circ p = p$.
A special example is the diagonal, denoted by $\Delta$.
Finally denote by $M(F)$ the category of {\bf effective Chow motives}, where
objects are pairs $(X,p)$ with $X \in V(F)$ and $p \in Mor(X,X)$ a projector.
The morphisms are described by $Mor((X,p),(Y,q)):=q \circ Mor(X,Y) \circ p$.\\

\begin{definition} Let $M=(X,p) \in M(F)$. Define
$$ CH^i(M):= p_*CH^i(X) \otimes \Q $$
\end{definition}

\begin{definition} 
Let $X \in V(F)$ be a smooth projective variety of dimension $d$.
We say that $X$ has a {\bf Murre decomposition}, if there exist projectors
$p_0,p_1,...,p_{2d}$ in $CH^d(X \times X) \otimes \Q$ such that the 
following properties hold (modulo rational equivalence for (1) and (2)):\\
(1) $p_j \circ p_i = \delta_{i,j} \cdot p_i $\\
(2) $\Delta = \sum p_i $\\
(3) In cohomology the $p_i$ induce the $(2d-i,i)-$th K\"unneth
component of the diagonal.\\
(4) $p_0,...,p_{j-1}$ and $p_{2j+1},..,p_{2d}$ act trivially on $CH^j(X) 
\otimes \Q$.\\
(5) If we put $F^0 CH^j(X) \otimes \Q = CH^j(X) \otimes \Q$ and inductively
$F^k CH^j(X) \otimes \Q := Ker(p_{2j+1-k} \mid_{F^{k-1}}) $, then this
descending filtration does not depend on the particular choice of the $p_i$.\\
(6) Always $F^1 CH^j(X) \otimes \Q = CH^j_{hom}(X) \otimes \Q$. 
\end{definition}

The motives $(X,p_i)$ are traditionally denoted by
$h^i(X)$ and we write $h(X)=h^0(X)+...+h^{2d}(X)$.  
In (6) one also wants to have that $F^2 CH^j(X) \otimes \Q$ is
the kernel of the cycle class map in rational Deligne cohomology, but this is
very hard to verify in general. \\
(1) - (6) have been proved for curves,
surfaces (\cite{6}), abelian varieties (\cite{7}) 
and certain varieties close to projective
varieties. Recently B. Gordon and 
J. Murre \cite{GM} computed the Chow motive of elliptic modular varieties 
using work of A. Scholl \cite{Tony}.\\
S. Saito has proposed a filtration in \cite{8} 
which has property (6).
Manin (\cite{3}) and Murre (\cite{6}) 
have quite generally defined
$p_0,p_1,p_{2d-1},p_{2d}$ for every $X$. A. Scholl has refined this
in \cite{Tony} to have also the property
that $p_1=p_{2d-1}^t$, where $p^t$ denotes a transpose of a projector $p$.
Murre has formulated the following \\
\ \\
{\bf Conjecture:} {\sl Every smooth projective F-variety $X$ admits 
a Murre decomposition.}\\

J. Murre has also studied the case of a product of a curve with a surface
where one in fact has a Murre decomposition. Inspired by this, we have
tried to construct projectors in the following situation: Let $f:Y \to S$
be a morphism from a smooth 3-fold $Y$  
to a smooth surface $S$ with connected fibers. 
Choose a smooth hyperplane section 
$i: Z \hookrightarrow Y$ and let $h=f|_Z$. 
Look the following cycles
$$\pi_{i0}:={1 \over m} (i \times 1)_* (h\times f)^*\pi _i(S),$$
$$\pi_{i2}:={1 \over m} (1 \times i)_* (f\times h)^*\pi _i(S),$$
in $CH^3(Y\times Y)\otimes\Q$. Here the $\pi_i(S)$ are  
orthogonal projectors of a 
Murre decomposition of $S$ as constructed by Murre (\cite{6}) 
and $m$ is the number of points on a general fiber of $h$. These cycles are 
not orthogonal in general but we are able to construct orthogonal projectors
$\pi_0,\cdots, \pi_6$ in the following way:

$$\pi_0:=\pi_{00},\quad \pi_1:=\pi_{10}, \quad 
\pi_2:=\pi_{20}+\pi_{02}-\pi_{20}\cdot\pi_{02}$$
$$\pi_4:=\pi_{40}+\pi_{22}-\pi_{40}\cdot\pi_{22}, \quad 
\pi_5:=\pi_{32}, \quad \pi_6:=\pi_{42}, \quad 
\pi_3:=\Delta-\sum_{i\ne 3}\pi_i
$$
The $\pi_j$ do not operate in the right way on cohomology, but if
all higher direct images sheaves $R^if_*\cO_Y$ vanish for $i \ge 1$,
they can be modified to form a Murre decomposition.
In particular a suitable blow up $Y$ of any smooth {\bf uniruled} 3-fold $X$
over a subfield of the complex numbers has this property.
Recall that a 3-fold $X$ is called uniruled, if there exists
a dominant rational map $\varphi: S \times {\PS}^1 - - - \to X$ for some
smooth projective surface $S$. By a theorem of Mori and 
Miyaoka (\cite{3}), this is equivalent to saying that $X$ 
has Kodaira dimension $-\infty$.
There is no structure theorem for these varieties which is as simple as in 
the case of ruled surfaces, but there is a version in the category of 
3-folds with ${\Bbb Q}$-factorial and terminal singularities (\cite{4}) 
stating that $X$ is birationally equivalent 
to a 3-fold $Y$ which has a fiber structure with rationally connected fibers
over a base variety which can be a point, a smooth curve or a
normal surface. Using this and suitable modifications of the
projectors above we can therefore prove:
\begin{thm} Let $X$ be a smooth uniruled complex projective 3-fold. 
Then $X$ admits a Murre decomposition.
\end{thm}
In the proof of this theorem, which makes heavy use of Fulton's
machinery of intersection theory, the Murre decomposition provides the 
following description of the {\bf Chow motive } of a complex 
uniruled 3-fold $X$ (ignoring torsion):\\
\[
\begin{tabular}{|l|c|c|c|c|c|c|c|c|c|} \hline
Motive $M$ &  $h^0(X)$  & \ $h^1(X)$ \ & \ $h^2(X)$ \ & \ $h^3(X)$ \ & \ 
$h^4(X)$ \ & \ $h^5(X)$ \ &  $h^6(X)$   \ \\ \hline
$CH^0(M)$ &  $\Z$  & \ 0 \ & \ 0 \ & \ 0 \ & \ 0 \ & \ 0 
\ &  $0$   \ \\ \hline
$CH^1(M)$ &  $0$  &  ${\rm Pic}^0(X)$  & \ ${\rm NS}(X)$ \ & \ $0$ \ 
& \ $0$ \ & \ $0$ \ & $0$   \ \\ \hline
$CH^2(M)$ &  $0$ & \ $0$ \ & \ ${\rm Ker}(\psi)$  \ & \ ${\rm Im}(\psi)$  \ 
& $H^{2,2}(X,\Z)$ & \ $0$ \ &  $0$   \ \\ \hline
$CH^3(M)$ &  $0$  & \ $0$ \ & \ $0$ \ & \ $0$ \ &  
${\rm Ker}(alb_X)$  
& \ $Alb(X)$ \ & $\Z$ \ \\ \hline
\end{tabular}
\]\\
\ \\
Here $\psi: CH^2_{\rm hom}(X) \to J^2(X)$ is the Abel-Jacobi map.
We hope that our approach may also be used to construct projectors in 
other situations. \\

{\bf Acknowledgements}: It is a pleasure to thank H. Esnault, 
B. Gordon, J. Murre and E. Viehweg for several discussions. 
The DFG and the universities of Essen and Leiden 
have supported the authors during this project.
\newpage    

\section{\bf Projectors for special varieties}
\bigskip 
The easiest case in which one has a Murre decomposition is the 
case of projective space, 
because there $H^{2k+1}(X,\C)=0$ for all $k\ge 0$ and the other groups
admit a basis represented by algebraic cycles. One has a more general theorem:\\

\begin{thm} \label{thm1}
Let $X$ be a smooth variety of dimension $n$ and assume that 
for certain $1\le q\le n-1$ there is a basis $\{E_1,\cdots ,E_t\}$ of
$H^{2q}(X,\Q)$ and a basis $\{\ell_1, \cdots ,\ell_t\}$
of $H^{2(n-q)}(X,\Q)$ represented by classes of algebraic
cycles. Then: \\
a) There exists a matrix $B=(b_{ij})\in {\bf GL_n}(\Q)$ such that the cycle
$p=\sum b_{ij}(\ell_i\times E_j)\in CH^n(X\times X)\otimes\Q$ operates as the
identity on $H^{2q}(X,\Q)$.\\
b) For the same choice of $b_{ij}$, 
$p^t=\sum b_{ij}(E_j\times \ell_i)
\in CH^n(X\times X)\otimes\Q$ operates as the
identity on $H^{2(n-q)}(X,\Q)$.\\
c) Both cycles, $p$ and $p^t$ are idempotent and therefore projectors.
\end{thm}

\dem Let $A=(E_i\cdot \ell_j)$ be the intersection matrix, then take
$B=A^{-1}$. \qed
Moreover, one can explicitely say how these projectors operate on 
cycles, namely:\\
\begin{prop} Let $p$ be as before and let $k\ne q$. Then, for all
$Z\in CH^k(X)\otimes\Q$ one has $p(Z)=0$ as an element of
$CH^k(X) \otimes \Q$.
\end{prop}
\dem By dimension reasons, as $p(Z)\in <E_i>\subset CH^q(X)\otimes\Q$. \qed
\begin{lemma} Let $p$ be as before and $Z\in CH^q(X)\otimes\Q$. If $[Z]$
denotes the homology class of $Z$ on $H^{2q}(X,\Q)$, then
$[p(Z)]=p([Z])=[Z]$.
\end{lemma}
\dem $p$ operates as the identity on $H^{2q}(X,\Q)$ and
$p(Z)=\sum b_{ij}(\ell_i\cdot Z)E_j$. \qed
\begin{cor} Let $p$ be as before, then
$({\rm Ker}\; p)\cap CH^q(X)\otimes\Q =CH_{\rm
hom}^q(X)\otimes\Q$.
\end{cor}
{\bf Examples:} Smooth Fano 3-folds and Calabi-Yau 3-folds have the property
that the Hodge numbers $h^{i,0}$ are always zero for $i=1,2$ and therefore
theorem \ref{thm1} applies. Another example is a del Pezzo fibration 
$f: X \to B$: 
to illustrate this, let $\ell$ be the extremal 
rational curve, $F$ a general fiber, 
$Y$ be a section of $|-mK_X|$, $C$ a twofold intersection in 
the linear system $|Y|$ and hence a 
multisection of $f$ over $B$, such that $C$ is a smooth 
curve dominating $B$. $H^2(X,\Q)$ is free of rank two. Then 
theorem \ref{thm1} produces the following projector
$$p_2:= {1 \over r}(C \times F)+{1 \over m }({\ell} \times Y)-
{d \over {m \cdot r}}({\ell} \times F) $$
where $d=Y^3$ and $r:=(C.F)$. Note that $(-K_X.\ell)=1$.
$p_2$ is unique as a cycle up to the choices of $Y,C,F$ and 
$\ell$.

\section{\bf Murre decompositions of birational conic bundles}
\bigskip
Let $f:Y\longrightarrow S$ be a morphism from a smooth projective 3-fold $Y$
to a smooth projective surface $S$, such that every fiber of $f$ is 
rationally
connected  and the general fiber of $f$ is isomorphic to ${\PS}^1$. 
Choose a smooth hyperplane section $i:Z\hookrightarrow Y$ such that 
$h:=f_{|Z}:Z\longrightarrow S$
is surjective and generically finite. Then define cycles
$$\pi_{i0}:={1 \over m} (i \times 1)_* (h\times f)^*\pi _i(S),$$
$$\pi_{i2}:={1 \over m}(1 \times i)_* (f\times h)^*\pi _i(S),$$
in $CH^3(Y\times Y)\otimes\Q$ for $0\le i\le 4$. 
Here the $\pi_i(S)$ are the orthogonal projectors of a 
Murre decomposition of $S$ as constructed by Murre (\cite{6}) 
(and improved by A. Scholl in \cite{Tony} to have also the property
that $\pi_1=\pi_3^t$)
and $m$ is the number of points on a general fiber of $h$. The following is
our {\bf key result} in some sense:
\begin{lemma}\quad  \\ 
a) $\pi_{i0}\circ \pi_{j0} = \delta_{ij}\pi_{i0}$ \\
b) $\pi_{i2}\circ \pi_{j2} = \delta_{ij}\pi_{i2}$ \\
c) $\pi_{j2}\circ \pi_{i0} = 0$ 
\end{lemma}

\dem a) Using the projection formula and the theory of Gysin maps for
l.c.i. morphisms from \cite[prop.6.6 (c)]{Fu} in the 
following diagram
$$
\matrix{Y \times Y  \times Y & \to & Y \times Y \cr
\downarrow && \downarrow \cr
Z \times Y \times Y & \to& Z \times Y \cr
\downarrow && \downarrow \cr
Z\times S \times Y &\to& Z \times Y \cr
\downarrow && \downarrow \cr
S \times S \times S & \to & S \times S }
$$
where the vertical maps are canonical l.c.i. morphisms, one obtains: \\
\begin{small}
$\pi_{i0}\circ \pi_{j0}$\\
$ ={1 \over m^2}(pr_{13}^{Y\times Y\times Y})_*
((i\times 1)_*((h\times f)^*(\pi_j(S))\times Y \cap
Y\times (i\times 1)_*((h\times f)^*(\pi_i(S))))$\\
$={1 \over m^2}(pr_{13}^{Y\times Y\times Y})_*
((i\times 1\times 1)_*(h\times f\times f)^*(\pi_j(S)\times S) \cap
(1\times i\times 1)_*(f\times h\times f)^*(S\times \pi_i(S)))$\\
$={1 \over m^2}(pr_{13}^{Y\times Y\times Y})_*
(i\times 1\times 1)_*[(h\times f\times f)^*(\pi_j(S)\times S) \cap
(i\times 1\times 1)^*(1\times i\times 1)_*(f\times h\times f)^*
(S\times \pi_i(S))]$\\
$={1 \over m^2}(i\times 1)_*(pr_{13}^{Z\times Y\times Y})_*
[(h\times f\times f)^*(\pi_j(S)\times S) \cap
(1\times i\times 1)_*(i\times 1\times 1)^*(f\times h\times f)^*
(S\times \pi_i(S))]$\\
$={1 \over m^2}(i\times 1)_*(pr_{13}^{Z\times Y\times Y})_*
[(h\times f\times f)^*(\pi_j(S)\times S) \cap
(1\times i\times 1)_*(h\times h\times f)^*(S\times \pi_i(S))]$\\
$={1 \over m^2}(i\times 1)_*(pr_{13}^{Z\times Y\times Y})_*(1\times i\times 1)_*
[(1\times i\times 1)^*(h\times f\times f)^*(\pi_j(S)\times S) \cap
(h\times h\times f)^*(S\times \pi_i(S))]$\\
$={1 \over m^2}(i\times 1)_*(pr_{13}^{Z\times S\times Y})_*(1\times h\times 1)_*
(h\times h\times f)^*[(\pi_j(S)\times S) \cap (S\times \pi_i(S))]$\\
$={1 \over m}(i\times 1)_*(pr_{13}^{Z\times S\times Y})_*
(h\times 1\times f)^*[(\pi_j(S)\times S) \cap (S\times \pi_i(S))]$ 
\hfill (\cite[prop. 6.6 (c)]{Fu}) \\
$={1 \over m}(i\times 1)_*(h\times f)^*(pr_{13}^{S\times S\times S})_*
(\pi_j(S)\times S) \cap (S\times \pi_i(S))$\\
$={1 \over m}(i\times 1)_*(h\times f)^*
(\pi_i(S) \circ \pi_j(S))=\delta_{ij}\pi_{i0}$.\\
\end{small}
Similarly one proves b).\\
c) As before, one finds that\\
\begin{small}
$\pi_{j2}\cdot \pi_{i0}$\\
$={1 \over m^2}(i\times i)_*(pr_{13}^{Z\times S\times Z})_*
(1\times f\times 1)_*
[(1\times f\times 1)^*(h\times 1\times h)^*
(\pi_i(S)\times S\cap S\times \pi_j(S))\cap (Z\times Y\times Z)]$\\
$={1 \over m^2}(i\times i)_*(pr_{13}^{Z\times S\times Z})_*
[(h\times 1\times h)^*(\pi_i(S)\times S\cap S\times \pi_j(S))\cap 
(1\times f\times 1)_*(Z\times Y\times Z)]= 0$\\
\end{small}
because $(1\times f\times 1)_*(Z\times Y\times Z)= 0$ due to
dimension reasons.  \qed
\ \\
Define now a set of cycles $\pi_0,\cdots,\pi_6$ in the following way:\\
\begin{small}
$\pi_0:=\pi_{00}$\\
$\pi_1:=\pi_{10}$\\
$\pi_2:=\pi_{20}+\pi_{02}-\pi_{20} \circ \pi_{02}$\\
$\pi_4:=\pi_{40}+\pi_{22}-\pi_{40} \circ \pi_{22}$\\
$\pi_5:=\pi_{32}$\\
$\pi_6:=\pi_{42}$\\
$\pi_3:=\Delta-\sum_{i\ne 3}\pi_i$
\end{small}
\begin{cor} 
The $\pi_j$ defined above form a set of orthogonal projectors 
such that $\pi_k=\pi_{6-k}^t$ .
\end{cor}

\begin{thm} \label{thm2} $\pi_i=\delta_{ij}$ on 
$\cases{f^*H^j(S,\Q)  \quad if \quad j = 0,1 \cr
        f^*H^j(S,\Q) \oplus \Q\cdot [Z]  \quad if \quad j=2 \cr
        f^*H^j(S,\Q) \oplus [Z]\cdot f^*H^2(S,\Q) \quad if \quad j=4 \cr
                            [Z]\cdot f^*H^3(S,\Q) \quad if \quad j=5 \cr
                            [Z]\cdot f^*H^4(S,\Q) \quad if \quad j=6}$
\end{thm}

\dem First note that one has the equation:
$\pi_{i0}(f^*\alpha)= 
{1 \over m}(i \times 1)_* (h\times f)^*\pi _i(S)(f^*\alpha)$\\
$={1 \over m}(pr^{Y \times Y}_2)_*[(i \times 1)_* 
(h\times f)^*\pi _i(S)\cap (f^*\alpha\times Y)]$\\
$ ={1 \over m}(pr^{Y \times Y}_2)_*(i \times 1)_* 
[(h\times f)^*\pi _i(S)\cap (i\times 1)^*(f^*\alpha\times Y)]$\\
$ ={1 \over m}(pr^{Y \times Y}_2)_*(i \times 1)_* 
(h\times f)^*[\pi _i(S)\cap \alpha\times S]$\\
$ ={1 \over m}(pr^{Z \times Y}_2)_* (h\times f)^*[\pi _i(S)
\cap \alpha\times S]$\\
$ ={1 \over m}(pr^{S \times Y}_2)_* (h\times 1)_*(h\times f)^*
[\pi _i(S) \cap \alpha \times S)]$\\
$ =(pr^{S \times Y}_2)_*(1\times f)^*[\pi _i(S)\cap \alpha\times S)]$\\
$ =f^*(pr^{S \times S}_2)_*[\pi _i(S)\cap \alpha\times S]
=f^*\pi_i(S)(\alpha)$.\\
Therefore $\pi_{i0}$ operates as $\delta_{ij}$ on $f^*H^j(S)$, 
proving the assertion for $\pi_0$ and $\pi_1$.\\
On the other hand, using projection formula, one gets\\
\noindent
$\pi_{i2}(f^*\alpha)={1 \over m}(pr^{Y\times Y}_2)_*
[ (1 \times i)_* (f\times h)^*\pi _i(S)\cap (f^* \alpha\times Y)$ \\
$\quad\quad\quad\quad = {1 \over m} i_* (pr^{S\times Z}_2)_*
(f\times 1)_*[(f\times 1)^*(1\times h)^*(\pi_i(S)\cap (\alpha\times S))
\cap (Y\times Z)]$\\
$\quad\quad\quad\quad = {1 \over m} i_* (pr^{S\times Z}_2)_*
[(1\times h)^*(\pi_i(S)\cap (\alpha\times S))
\cap (f\times 1)_*(Y\times Z)]=0$,\\
since $(f\times 1)_*(Y\times Z)=0.$\\
Take any $D\in H^k(S,\Q)$ with $k=0,2,3,4$ and consider 
$C:=i_*h^*(D)$. Observe that 
$[C]=f^*(D)\cdot [Z]$. The same computation as above in
cohomology shows that 
$$\pi_{i2}([C])=:{1 \over m}(pr_2^{Y\times Y})_*
[(1\times i)_*(f\times h)^*\pi_i(S)\cap [C]\times [Y]]=
i_*h^*(\pi_i(S)(D))$$
As the $\pi_i(S)$ induce the K\"unneth decomposition of 
$\Delta_S$ on cohomology,
it follows that $\pi_i(S)([D])=\delta_{ik}([D])$ and therefore one gets
$\pi_{i2}([C])=\delta_{ik}[C]$.\\
Moreover, a similar argument together with Chow's moving lemma shows that\\
$\pi_{i0}([C])={1 \over m}(pr_2^{Y\times Y})_*
[(i \times 1)_*(h \times f)^*\pi_i(S)\cap [C] \times [Y]]$\\
$={1 \over m}(pr_2^{Y\times Y})_*
(i \times 1)_*[(h \times f)^*\pi_i(S)\cap (i \times 1)^* [C] \times [Y]]$\\
$={1 \over m}(pr_2^{Z \times Y})_*
[(h \times f)^*\pi_i(S)\cap [C \cap Z] \times [Y]]$\\
$={1 \over m}(pr_2^{S \times Y})_*(h \times 1)_* 
[(h \times 1)^*(1 \times f)^*
\pi_i(S)\cap [C \cap Z] \times [Y]]$\\
$={1 \over m}(pr_2^{S \times Y})_*[(1 \times f)^*
\pi_i(S)\cap h_* [C \cap Z] \times [Y]]$\\
$={1 \over m} f^* (pr_2^{S \times S})_*[ \pi_i(S) \cap h_* 
[C \cap Z] \times [S]]$\\
$={1 \over m} f^* \pi_i(S)(h_* [C \cap Z])=0$,\\
if $i \ne k+2$. As a consequence one also gets
$\pi_{i0} \circ \pi_{j2}([C])=\delta_{jk}\pi_{i0}([C])$, which proves the
assertion for $\pi_2,\pi_4,\pi_5$ and $\pi_6$ and the theorem. \qed
\ \\
Now assume additionally that $f: Y \to S$ is a desingularization of a conic
bundle morphism $f': X' \to S' $ in the sense of \cite{4}, 
i.e. there is a commutative diagram
$$\matrix{ Y & {\buildrel f \over \longrightarrow} & S \cr
\quad \downarrow \sigma && \quad \downarrow \tau  \cr
X' & {\buildrel f' \over \lra} & S' }
$$
with blow-up morphisms $\sigma, \tau$. Also we assume $Z \subset Y$ is
a sufficiently general smooth hyperplane section 
of $Y$ that dominates $S$. 

Then we can choose irreducible
divisors $H_1,...,H_r$ in $Y$ such that $H_1=Z$ and 
$$ H^{1,1}(Y,\Q)= \bigoplus_{i=1}^r \Q[H_i] 
$$
form a basis of $H^{1,1}(Y,\Q)$ and such that 
$f_*H_i=0$ in $CH^0(S)$ for $i \ge 2$, i.e. $H_i$ is exceptional with
respect to $f$ for $i \ge 2$. 

\begin{lemma} \label{le1}
For every cycle $W$ one has $\pi_{20}(W)={1 \over m} 
f^* \pi_2(S) (h_* (W \cap Z)) \in f^* CH^*(S) \otimes \Q$.
Let $W$ be a cycle with $f_*(W) =0$. 
Then $\pi_{02}(W)=0$ already in the Chow group of $Y$. 
\end{lemma} 
\dem 
$\pi_{02}(W)={1 \over m}(pr_2^{Y\times Y})_*
[(1\times i)_*(f\times h)^*\pi_0(S)\cap (W \times Y) ]$\\
$={1 \over m} i_* (pr_2^{S \times Z})_*[(1 \times h )^* 
\pi_0(S) \cap (f \times 1)_*(W \times Z)] =0 $ \\
by \cite[prop.6.6 (c)]{Fu} and  
since $f_*(W)=0 \in CH^*(S)$. \\
On the other hand\\
$\pi_{20}(W)={1 \over m}(pr_2^{Y\times Y})_*
[(i \times 1)_* (h \times f)^*\pi_2(S) \cap (W \times Y) ]$\\
$={1 \over m}(pr_2^{Z \times Y})_*
[(h \times f)^*\pi_2(S) \cap ((W \cap Z) \times Y) ]$\\
$={1 \over m} (pr_2^{S \times Y})_*[(1 \times f)^* 
\pi_2(S) \cap (h \times 1)_*((W \cap Z) \times Y)] $\\
$={1 \over m} (pr_2^{S \times Y})_*(1 \times f)^*[ 
\pi_2(S) \cap h_*(W \cap Z) \times S)] $\\
$={1 \over m} f^* (pr_2^{S \times S})_*
[\pi_2(S) \cap h_*(W \cap Z) \times S ]  $\\
$={1 \over m} f^* \pi_2(S)(h_*(W \cap Z)) \in f^* CH^*(S) \otimes \Q$.
\qed

\begin{cor} $\pi_2(Y)(H_i)={1 \over m} f^*(h_* (H_i \cap Z)) 
\in f^* CH^1(S) \otimes \Q $ for $i \ge 2$.
\end{cor}

By theorem \ref{thm2} $\pi_2(Y)$ operates as zero on ${\rm Pic}^0(Y)$, 
therefore the image of $\pi_2(Y)$ in $CH^1(Y) \otimes \Q$ is a finite
dimensional vector space.
By changing our generators $H_i$ above modulo classes in ${\rm Pic}^0(Y)=
f^*{\rm Pic}^0(S)$, we may assume that they  
generate ${\rm Im}(\pi_2) \subset CH^1(Y) \otimes \Q$. 
Then we write uniquely 
$$\pi_2(Y)(H_i)=\sum_k a_{i,k} H_k \in CH^1(Y) \otimes \Q$$ 
with a matrix $A=(a_{i,k}) \in {\rm Mat}(r \times r,\Q)$. 
$\pi_2(Y)$ being a projector implies that $A^2=A$. 
Choose algebraic cycles $\ell_1,...,\ell_r $ such that $\ell_1=F$, a general
fiber of $f$, and such that
their cohomology classes form a basis of $H^{2,2}(Y,\Q)$. 
By Poincar\'e duality the intersection matrix
$M=(m_{i,j}):=
(\ell_1,...,\ell_r)^{T}(H_1,...,H_r)$
has nonzero determinant.\\
We define 
$$ q_2:= \pi_2(Y) + \sum b_{i,j}(\ell_i \times H_j) - 
\sum b_{i,j}(\ell_i \times H_j) \circ \pi_2
$$
with some matrix $B=(b_{i,j}) \in {\rm Mat}(r \times r,\Q)$.\\
 
\begin{lemma}{\label{le2}} If $B=M^{-1}({\bf 1} -A)$, then $q_2$ 
is a projector and operates as the
identity on $H^2(Y,\Q)$. 
\end{lemma}
\dem 
$\pi_2$ acts as the identity on $f^*H^2(S,\Q)$ by
theorem \ref{thm2}. The higher direct images $R^if_*\cO_Y$ 
vanish for $i \ge 1$ by \cite{4}. Therefore by the Leray spectral sequence
$H^2(Y,\cO_Y)=f^*H^2(S,\cO_S)$ and it is
enough to show that $q_2$ operates as the identity on $H^{1,1}(Y,\Q)$ too.
But $q_2$ acts via the matrix $MB + A +BA$
on $H^{1,1}(Y,\Q)$ with respect to 
the basis $\{H_i\}$. Now $\pi_2^2=\pi_2$ and we get
$A^2=A$ and therefore $BA=0$. By definition of $B$, we obtain
that $MB + A +BA= M(M^{-1}({\bf 1}-A)) +A={\bf 1}$. \\
To show that $q_2$ is a projector, let us write 
$q_2=\pi_2+ \beta - \beta \pi_2$. Note that $\beta \beta = \beta$,
since $BMB=B$. From $BA=0$ we deduce that $\pi_2 \beta=0$.
Therefore\\
$q_2 \circ q_2= \pi_2^2 +\beta^2 + \beta \pi_2 \beta \pi_2 + \pi_2 \beta 
- \pi_2 \beta \pi_2 + \beta \pi_2 - \beta \beta \pi_2 - \beta \pi_2 \pi_2
-\beta \pi_2 \beta = \pi_2 +\beta -\beta \pi_2 = q_2$ is a projector.   
\qed 

\begin{thm} \label{thm3} The following cycles
$ p_0(Y):=\pi_0(Y), \ p_1(Y):=\pi_1(Y) $, \\ 
$ p_2(Y):= q_2 -\pi_1(Y) \circ \sum b_{i,j}(\ell_i \times H_j) 
-\pi_1(Y) \circ \sum b_{i,j}(\ell_i \times H_j) \circ \pi_2(Y) $\\
$ p_4:= p_2^{tr}(Y), \quad p_5(Y):= \pi_5(Y),
\quad p_6(Y):=\pi_6(Y), \quad p_3(Y):=\Delta - \sum_{i \ne 3} p_i $\\
define orthogonal projectors, which satisfy properties (1)-(6) of a 
Murre decomposition. 
\end{thm}

\dem By lemma \ref{le2} above, (1),(2) and (3) are straightforward. \\
To prove (4),(5) and (6) 
for $j=1$, note that ${\rm Pic}(Y) \otimes \Q =f^*{\rm Pic}^0(S) \otimes \Q 
\oplus \bigoplus_i \Q \cdot H_i$. By theorem \ref{thm2} above,
$p_1$ operates on $Pic^0(Y) \otimes \Q=f^*Pic^0(S) \otimes \Q$ 
as the identity and trivially on
$\bigoplus_i \Q \cdot H_i$. Vice versa $p_2$ is the identity on
$\bigoplus_i \Q \cdot H_i$ and zero on $f^*Pic^0(S) \otimes \Q$, because it
acts trivially on $f^*H^1(S,\Q)$. 
All the other projectors are zero on $CH^1(Y) \otimes \Q$.  
Therefore we get (4)-(6) for $j=1$ with $F^2CH^1(Y) \otimes \Q=0$.\\
For $j=2$, property (4) follows from the analogous assertion for $S$. By
construction $F^1 CH^2(Y) \otimes \Q={\rm Ker}(p_4)=
CH^2_{\rm hom}(Y) \otimes \Q$. Then
$F^2CH^2(Y) \otimes \Q={\rm Ker}(p_3)\cap{\rm Ker}(p_4)={\rm Im}(p_2)=
{\rm Im}(\pi_2(Y))$.\\
Now we show that $F^2CH^2(Y) \otimes \Q \cong f^* F^2 CH^2(S) \otimes \Q$: 
$\pi_{02}$ operates as zero on $CH^2(Y)$ 
by Chow's moving lemma and if $C$ is any
curve homologous to zero on $Y$, then by Lemma \ref{le1}, 
$\pi_{20}(C)= f^* h_*(C \cap Z) \in 
f^* F^2 CH^2(S) \otimes \Q$.\\
This proves that
$F^2CH^2(Y) \otimes \Q \subset f^* F^2CH^2(S) \otimes \Q$, but
since $\pi_2(Y)$ operates as the identity on every fiber of $f$, we get
equality. This is then independent of all choices, because this is the case
for $F^2 CH^2(S)$ by \cite{6}. 
Finally $F^3 CH^2(Y) \otimes \Q=0$, since $p_2$ acts as the
identity on $F^2CH^2(Y) \otimes \Q ={\rm Im}(p_2)$. 
Hence we get (5) and (6) for $j=2$.\\
Finally consider $CH^3(Y)$: Clearly $F^1CH^3(Y) \otimes \Q=Ker(\pi_6)=
CH^3_{\rm hom}(Y) \otimes \Q$. Further 
$F^2CH^3(Y) \otimes \Q = Ker (\pi_5|_{F^1CH^3(Y) \otimes \Q})$ and we 
claim that $F^2 CH^3(Y) \otimes \Q \cong Ker(alb_Y) \otimes \Q$, where
$alb_Y: CH^3(Y)_{\rm hom} \to Alb(Y)$ is the Albanese map. But
there is a commutative diagram
$$\matrix{ CH^3(Y)_{\rm hom} & \to & Alb(Y) \cr
f_* \downarrow && \downarrow f_* \cr
CH^2(S)_{\rm hom} & \to & Alb(S) }
$$
Both vertical maps are isomorphisms. To compute
$F^2CH^3(Y) \otimes \Q $ we take any closed point $P$ in $Y$ and compute that
$f_* \pi_5(P)= f_* {1 \over m} i_* h^* (\pi_3(S)(P))= \pi_3(S)(f_*(P))$.\\
This shows that $f_* F^2CH^3(Y) \otimes \Q  \cong 
F^2 CH^2(S) \otimes \Q \cong Ker(alb_S) \otimes \Q $ by \cite{6}. Therefore
$F^2CH^3(Y) \otimes \Q \cong Ker(alb_Y) \otimes \Q$, which is 
independent of all choices again by \cite{6}. Finally 
$F^3CH^3(Y) \otimes \Q = 0$, since if $P=\sum a_i P_i$ is a zero
cycle on $Y$ with $\sum a_i=0$, then 
$f_* \pi_4(P)=f_* \pi_{20}^{t}(P) + f_* \pi_{02}^{t}(P) =
f_*{1 \over m} (1 \times i)_*(f \times h)^*\pi_2(S) (P) +
f_*{1 \over m} (i \times 1)_* (h \times f)^* \pi_4(S) (P)$. 
But $\pi_4(S)=S \times e$ , hence the last term is zero  and the first 
term becomes $\pi_2(S)(f_*P)$. But $\pi_2(S)$ acts as the identity
on $F^2 CH^2(S) \otimes \Q$. Thus $f_* F^3 CH^3(Y) \otimes \Q
\subset F^3 CH^2(S) \otimes \Q = 0$. \\
This finishes the proof of the theorem. \qed

\newpage
\section{Murre decompositions of uniruled 3-folds}

Let $k=\C$. By a 3-fold we just mean a normal 3-dimensional complex 
variety.\\
\begin{definition} A 3-fold $X$ is called {\bf uniruled}, if there exists
a dominant rational map $\varphi: S \times {\PS}^1 - - - \to X$ for some
surface $S$.
\end{definition}

\begin{thm} (\cite{3}): A smooth projective 3-fold $X$ is uniruled
if and only if it has Kodaira dimension $- \infty$, i.e. no multiple of
$K_X$ has sections.
\end{thm}

\begin{thm} (\cite{4}): Let $X$ be a uniruled 3-fold
with only $\Q$-factorial terminal singularities. Then there exists a 
birational mapping $r : X ---\to Y$ which is a composition of flips and 
divisorial contractions, such that $Y$ has an extremal ray $R$ whose 
extremal contraction map $f: Y \to Z$ satisfies one of the following cases:\\
(a) dim(Z)=0, $Y$ is a $\Q$-Fano 3-fold with $\rho(Y)=1$, i.e. $-mK_Y$ is
an ample Cartier divisor for some $m \ge 1$ and the divisor class group is
free with one generator.\\
(b) $Z$ is a smooth curve and $Y$ is a del Pezzo fibration over $Z$, i.e.
the general fibre of $f$ is a del Pezzo surface.\\
(c) $Z$ is a surface with at most quotient singularities and $Y$ is a 
conic bundle over $Z$.\\
In cases (b) and (c) the reduced preimage of any irreducible divisor is again
irreducible. 
\end{thm}

\begin{thm} Let $X$ be a smooth complex 
uniruled 3-fold. Then $X$ admits a Murre decomposition.
\end{thm}
\bigskip
\dem Since $X$ is uniruled, it is birational to one of the following 
varieties:\\  
(a) A $\Q$-Fano 3-fold $Y$ with $\rho(Y)=1$, i.e. $-mK_Y$ is
an ample Cartier divisor for some $m \ge 1$ and the divisor class group is
free with one generator.\\
(b) A del Pezzo fibration over a smooth curve. \\
(c) A conic bundle over a normal 
surface with at most quotient singularities.\\
In cases (a), (b) $H^2(X,\Q)$ and $H^4(X,\Q)$ are 
generated by classes of 
algebraic cycles. Thus we define $p_0(X)=\{e\}\times X$
and $p_6(X)=X \times\{e\}$ for some rational point $e\in X$, $p_1(X)$ and
$p_5(X)$ as in \cite{6} and $p_2(X)$ and $p_4(X)=p_2(X)^{tr}$ as in 
theorem \ref{thm1}. Then it is immediate to verify all properties 
(2)-(6) similar to the proof of \ref{thm3} while property (1) can be achieved 
like in \cite[remark 6.5.]{6}, by the non-commutative Gram-Schmidt process.\\
In case (c) we may assume that after blowing up $X$ along several
smooth subvarieties, there is a situation as in the previous section:\\
Let $\varphi: Y \to X$ be the blow-up and assume that $f:Y \to S$
is a morphism to a smooth surface $S$ with rationally connected fibers.
Take the projectors $p_0(Y),...,p_6(Y)$ as defined in the last section.\\
To define the projectors for $X$, 
consider the graph $\Gamma_\varphi \subset
Y \times X$ of $\varphi$. Define 
$$p_i(X):=\Gamma_\varphi \circ p_i(Y) \circ \Gamma_\varphi^{tr}= 
(\varphi\times\varphi)_*(p_i(Y))$$
(by Liebermann's lemma \cite{2})
for $0 \le i \le 2$. We claim that all $p_i(X)$ are orthogonal projectors.\\
By induction on the number of blow-ups we may assume that there is 
just one blow-up along a smooth subvariety $W\subset X$.\\
Consider the canonical diagram 
$$\matrix{Y \times Y \times Y & {\buildrel pr_{13} \over \to} 
& Y \times Y \cr
\downarrow && 
\downarrow \cr
X \times Y \times X & {\buildrel pr_{13} \over \to} & X \times X}
$$
where the vertical maps are $\varphi \times 1 \times \varphi$ and
$\varphi \times \varphi$. 
Let $E$ be the exceptional divisor. Then we compute for $0 \le i,j \le 2$:\\
$p_i(X)\circ p_j(X)= (pr_{13})_*((\varphi\times id)_*p_j(Y)\times X\cap
X\times (id\times \varphi)_*p_i(Y))=$\\
$= (\varphi\times\varphi)_*
(pr_{13})_*(p_j(Y)\times Y\cap
Y\times (id\times \varphi)^*(id\times \varphi)_*p_i(Y))$\\
$= (\varphi\times\varphi)_*
(pr_{13})_*(p_j(Y)\times Y\cap
Y\times (p_i(Y)+(id\times j)_*Q_{i,j}))$ \\
$=(\varphi\times\varphi)_* (pr_{13})_*(p_j(Y)\times Y\cap Y\times p_i(Y))+
(\varphi\times\varphi)_* 
(pr_{13})_*(p_j(Y)\times Y\cap Y \times (id\times j)_*Q_{i,j})$\\
$=(\varphi\times\varphi)_* (p_i(Y) \circ p_j(Y)+
(pr_{13})_*(p_j(Y)\times Y \cap Y \times (id\times j)_*Q_{i,j})$\\
where $Q_{i,j}\in CH_3(Y\times E)$ and $j: E \hookrightarrow Y$ is
the inclusion. Hence \\
$\cC_i:=p_i(X)\circ p_i(X)-p_i(X)= (\varphi\times id)_*
(pr_{13})_*(p_i(Y)\times X \cap
Y \times (id\times i)_*(id\times \varphi^E)_*Q_{i,i}))$.\\
$p_i(Y)={1 \over m} (i \times 1)_*(h \times f)^*\pi_i(S) + T_i$
with $T_0,T_1=0$ and 
$T_2=\sum c_{ij}(\ell_i \times H_j)- 
\sum b_{i,j}(\ell_i \times H_j) \circ \pi_2(Y)$ for some integers $c_{i,j},
b_{i,j}$ which
is supported on $(Z \times Y) \cup (\ell_i \times Y)$. Therefore
$\cC_i$ is supported on $\varphi(Z) \times W$.
Here $i:W\to X$ is the inclusion and $\varphi^E: E \to W$ is the restriction
of $\varphi$ to $E$.\\
If $W$ is a point, $\cC_i=0$ by dimension reasons. If $W$ is a curve,
$\cC_i=a(\varphi(Z)\times W)$ with $a \in \Z$. But 
$\cC_i=p_i(X)\circ p_i(X)-p_i(X)$ 
operates as zero on the cohomology classes of every curve $T \in CH^2(X)$, 
since by Chow's moving lemma we
can choose $T$ to be disjoint from $W$ and use that
$p_i(Y)(T)=0$ in cohomology for $i=0,1,2$. 
Therefore $a=0$ and $p_i(X)$ is a projector.\\
For $ i \neq j$, 
$p_i(X)\circ p_j(X)=
(\varphi\times\varphi)_* 
(pr_{13})_*(p_j(Y)\times Y \cap Y \times (id\times j)_*Q_{i,j})$\\
since $p_i(Y)$ and $p_j(Y)$ are orthogonal. As above this implies that
$ p_i(X) \circ p_j(X)$ is supported on
$\varphi(Z) \times W$ for all $j$. By the same argument 
with Chow's moving lemma for $CH^2(X)$ as before, $p_i(X)\circ p_j(X)=0$. \\ 
Now define  
$$p_4(X)=p_2(X)^{tr}, p_5(X)=p_1(X)^{tr}, p_6(X)=p_0^{tr} \quad {\rm and} 
\quad p_3(X)=\Delta-\sum_{i \ne 3} p_i(X)$$ 
Properties (3)-(6) follow from theorem
\ref{thm3} together with the split exact sequences (\cite[prop. 6.7]{Fu})
$$ 0 \to CH_k(W) \to CH_k(E) \oplus CH_k(X) \to CH_k(Y) \to 0 $$
(1) and (2) can be obtained again via the 
Gram-Schmidt process. \qed

{\bf Mail Address:} 
{\tt Universit\"at Essen, Fachbereich 6, 45117 Essen, Germany \\
e-mail: mueller-stach@uni-essen.de, plar@xanum.uam.mx}

\end{document}